\documentclass[11pt, letterpaper]{amsart}
\usepackage{mathbbol}
\usepackage{amsmath, amsthm, amssymb}
\usepackage{txfonts}

\newtheorem{thm}{Theorem}[section]
\newcommand{\bt}{\begin{thm}}
\newcommand{\et}{\end{thm}}

\newtheorem{cor}[thm]{Corollary}
\newcommand{\bc}{\begin{cor}}
\newcommand{\ec}{\end{cor}}

\newtheorem{lem}[thm]{Lemma}
\newcommand{\bl}{\begin{lem}}
\newcommand{\el}{\end{lem}}

\newtheorem{prop}[thm]{Proposition}
\newcommand{\bp}{\begin{prop}}
\newcommand{\ep}{\end{prop}}

\newtheorem*{prop*}{Proposition}

\newtheorem{defn}[thm]{Definition}
\newcommand{\bd}{\begin{defn}}      
\newcommand{\ed}{\end{defn}}

\newtheorem{rmrk}[thm]{Remark}
\newcommand{\br}{\begin{rmrk}}
\newcommand{\er}{\end{rmrk}}

\newcommand{\thmref}[1]{Theorem~\ref{#1}}
\newcommand{\secref}[1]{Section~\ref{#1}}
\newcommand{\lemref}[1]{Lemma~\ref{#1}}

\newcommand{\corref}[1]{Corollary~\ref{#1}}
\newcommand{\propref}[1]{Proposition~\ref{#1}}
\newcommand{\remref}[1]{Remark~\ref{#1}}

\newcommand{\N}{\mathbb{N}}

\newcommand{\R}{\mathbb{R}}

\newcommand{\dist}{\operatorname{dist}}
\newcommand{\diam}{\operatorname{diam}}

\newcommand{\lm}{{\mathcal L}}
\newcommand{\md}{\operatorname{md}}

\newcommand{\image}{\operatorname{Im}}

\newcommand{\ohne}{\backslash}
\newcommand{\bdry}{\partial\hspace{-0.05cm}}

\newcommand{\lip}{\operatorname{Lip}}

\begin{document}

\title{Characterizations of metric trees and Gromov hyperbolic spaces}

\author{Stefan Wenger}

\address
  {Courant Institute of Mathematical Sciences\\
   251 Mercer Street\\
   New York, NY 10012}
\email{wenger@cims.nyu.edu}

\date{August 25, 2007}


\begin{abstract} 
In this note we give new characterizations of metric trees and Gromov hyperbolic spaces and generalize recent results of Chatterji and Niblo.
Our results have a purely metric character, however, their proofs involve two classical tools from analysis: Stokes' formula in  $\R^2$ and a Rademacher type differentiation theorem for Lipschitz maps. This analytic approach can be used to give characterizations
of Gromov hyperbolicity via isoperimetric inequalities with optimal constants.
\end{abstract}

\maketitle

\bigskip

\section{Statement of the main results}
A geodesic metric space is said to be a metric tree (or an $\R$-tree) if it is $0$-hyperbolic in the sense of Gromov, or in other words, if all its geodesic triangles are isometric to tripods. Metric trees are used in geometry, metric topology, geometric group theory and also in the geometry of Banach spaces (see e.g.~\cite{Preiss-Schechtman} for the latter).
They have the following property which can easily be verified (see \lemref{Lemma:hyp-distortion}): The intersection $B_1\cap B_2$ of any two closed balls $B_1, B_2$ is a ball or the empty set. 
In this note we prove that this property already characterizes metric trees and that, maybe somewhat surprisingly, already a weaker property for a geodesic metric space implies 
that it is a metric tree. Namely, we show: 
\bt\label{theorem:main-thm-trees}
  Let $(X,d)$ be a geodesic metric space and $\lambda\geq 1$. If for any two balls $B_1,B_2\subset X$ with non-empty intersection
  there exist $z,z'\in X$ and $\nu\geq 0$ such that
  \begin{equation*}
   B(z,\nu)\subset B_1\cap B_2\subset B(z',\lambda\nu) 
  \end{equation*}
  then $X$ is a metric tree.
\et
Here, $B(z,\nu):= \{x\in X: d(x,z)\leq \nu\}$. As a special case we obtain the following characterization of metric trees, initially conjectured by Chatterji-Niblo and, 
shortly after a preliminary version of our paper had appeared,  established independently in \cite{Chatterji-Niblo} with metric methods:

\bc
 Let $(X,d)$ be a geodesic metric space. Then $X$ is a metric tree if and only if the intersection of any two closed balls is a ball or the empty set.
\ec

As another application of  \thmref{theorem:main-thm-trees} we obtain a generalization of the main result of \cite{Chatterji-Niblo}.
\bc\label{corollary:gromov-charact}
 Let $(X,d)$ be a geodesic metric space and suppose there exist $\lambda, \delta>0$ 
 with the following property: For any two balls $B_1, B_2\subset X$ with non-empty intersection there exist $z, z'\in X$ and $r\geq 0$ such that
  \begin{equation}\label{equation:cor-distortion}
   B(z, r) \subset B_1\cap B_2 \subset B(z', \lambda r+\delta).
  \end{equation}
 Then $X$ is Gromov hyperbolic.
\ec
Conversely, if $X$ is a geodesic Gromov hyperbolic space then it is not difficult to show that \eqref{equation:cor-distortion} holds with $\lambda=1$ and some 
$\delta$ only depending on the constant of Gromov hyperbolicity. For a proof of this see \lemref{Lemma:hyp-distortion} or \cite{Chatterji-Niblo}. For the definition of Gromov
hyperbolicity see \secref{section:preliminaries}.

We note that in general,  \thmref{theorem:main-thm-trees} is false if $X$ is not geodesic. Indeed, if $X:= (\R, |x-y|^{\frac{1}{2}})$ then the intersection of any two closed
 balls in $X$ is either a ball or the empty set. However, there does not exist a bilipschitz embedding of $X$ into a metric tree.

A remarkable feature of our note is the analytic flavor of our proof of the purely metric statement of  \thmref{theorem:main-thm-trees}. 
This analytic approach has some advantages over a metric one and can for example be used
to give a characterization of geodesic Gromov hyperbolic spaces in terms of  a quadratic isoperimetric inequality with the sharp isoperimetric constant, see \cite{Wenger-sharp-constant}.

\subsection{Outline of the proof of \thmref{theorem:main-thm-trees}}
The starting point of the proof of our main theorem is the following simple but useful observation which provides yet another new characterization of metric trees.

\begin{prop*}
 A geodesic metric space $X$ is a metric tree if and only if for every Lipschitz loop $\gamma:S^1\to X$ and all Lipschitz functions $f, \pi:X\to \R$ 
 \begin{equation}\label{equation:trivloops-intro}
  \int_{S^1} (f\circ\gamma)(t)\cdot(\pi\circ\gamma)'(t)\,dt= 0.
 \end{equation}
\end{prop*}

This proposition will be proved in \secref{section:charact-trees} and will be used as follows: 
Let $X$ satisfy the hypotheses of \thmref{theorem:main-thm-trees}.
In \lemref{lemma:bicombing-property} and \corref{corollary:lip-connected} we first show that $X$ is Lipschitz $1$-connected, i.e., 
that for
every Lipschitz loop
$\gamma: S^1\to X$ there exists a Lipschitz extension $\tilde{\gamma}:B^2\to X$ to the unit disc $B^2\subset\R^2$. Next, we prove in \lemref{lemma:bilipschitz-piece} that
\begin{equation}\label{equation:det=0}
 \det\left(D_z[(f,\pi)\circ\tilde{\gamma}]\right)=0\quad\text{for almost every $z\in B^2$,}
\end{equation}
for all Lipschitz functions $f,\pi: X\to\R$. Here, $D_z[(f,\pi)\circ\tilde{\gamma}]$ denotes the classical derivative at $z$, which exists almost everywhere by Rademacher's theorem.
Roughly, \eqref{equation:det=0} holds because otherwise, $X$ would contain almost-isometric copies of pieces of a $2$-dimensional normed space by a generalized Rademacher 
theorem. This however is 
easily seen to contradict the assumption on the intersection of balls.
Finally, Stokes' theorem together with \eqref{equation:det=0} yields \eqref{equation:trivloops-intro}, implying that $X$ is a metric tree.
The generalized Rademacher theorem alluded to above is due to Kirchheim and independently to Korevaar-Schoen. In the appendix we will give a short and partially new proof 
of this theorem.

\bigskip

{\bf Acknowledgments:}
I would like to thank Indira Chatterji for a discussion and Bruce Kleiner for pointing out to me the possibility of proving \corref{corollary:gromov-charact} using our main theorem.
I am furthermore indebted to Theo B\"uhler for many suggestions and comments.

\section{Preliminaries}\label{section:preliminaries}
\subsection{Geodesic metric spaces, metric trees and Gromov hyperbolicity}
Let $(X,d)$ be a metric space and $a>0$. A map $c:[0,a]\to X$ is called a (constant speed) geodesic path from $c(0)$ to $c(a)$ if there exists $\lambda>0$ such that
$d(c(s), c(t))=\lambda|s-t|$ for all $s,t\in [0,a]$. If $\lambda=1$ then $c$ is called a geodesic path parameterized by arc-length. The metric space $X$ is called 
geodesic if any two points can be joined by a geodesic path. 

A geodesic triangle in $X$ consists of three points (the vertices) in $X$ and a choice of three geodesics (the edges) joining these points. 
Let $\Delta$ be a geodesic triangle and $\delta\geq0$. Denote by $x_1,x_2,x_3\in X$ its vertices and by $c_{k\ell}$ its edges (joining $x_k$ to $x_\ell$), which we may assume to be
parameterized by arc-length, $1\leq k<\ell\leq 3$.
There exist unique numbers $a_1,a_2,a_3\geq 0$ with the property that $d(x_k,x_\ell)= a_k+a_\ell$.
If $1\leq k<\ell\leq 3$ we define $c_{\ell k}:[0, a_k+a_\ell]\to X$ by $c_{\ell k}(t):= c_{k\ell}(a_k+a_\ell-t)$.
Then $\Delta$ is said to be $\delta$-thin if for every permutation $\sigma$ of $\{1, 2, 3\}$ 
 \begin{equation*}
  d(c_{\sigma(1)\sigma(2)}(t), c_{\sigma(1)\sigma(3)}(t))\leq \delta
 \end{equation*}
 for all $t\in [0, a_{\sigma(1)}]$.

\bd 
 A geodesic metric space $(X,d)$ is called Gromov hyperbolic if there exists $\delta\geq 0$ such that every geodesic triangle in $X$ is $\delta$-thin. If $\delta=0$ 
 then $X$ is called a metric tree.
\ed

Gromov hyperbolic spaces were introduced and first studied by Gromov in \cite{Gromov-hyperbolic} in the context of geometric group theory. See  for example
\cite{Bridson-Haefliger} or \cite{Ghys-delaHarpe} for a treatment of the topic.

\subsection{Lipschitz maps and metric derivatives}
The following notion of differentiability for metric space valued Lipschitz maps was introduced by Kirchheim in \cite{Kirchheim}.
Let $\varphi: U\to X$ be a Lipschitz map, where $U\subset\R^k$ is open and $(X,d)$ a metric space. 
The metric directional derivative of $\varphi$ at $z\in U$ in direction $v\in\R^k$ is defined by
\begin{equation*}
 \md\varphi_z(v):= \lim_{r\searrow 0}\frac{d(\varphi(z+rv),\varphi(z))}{r}
\end{equation*}
if this limit exists. 
We will make use of the following trivial remark in the proof of our main theorem.
\br\label{remark:det-degenerate}{\rm
 Let $\varphi: U\subset\R^k\to X$ and $\varrho: X\to \R^k$ be Lipschitz maps and $z\in U$. If $\varrho\circ\varphi$ is differentiable at $z$ and $\md\varphi_z(v)$ exists
 for all $v\in\R^k$ and is degenerate (that is $\md\varphi_z(v)=0$ for some $v\not=0$), then
 \begin{equation*}
  \det\left(D_z[\varrho\circ\varphi]\right)=0.
 \end{equation*}
}
\er

The following theorem was proved by Kirchheim in \cite{Kirchheim} and a similar statement by Korevaar-Schoen in \cite{Korevaar-Schoen}. We will provide a partially new proof
in the appendix.
\bt\label{theorem:metric-deriv-strong}
 Let $(X,d)$ be a metric space and $\varphi:U\to X$ a Lipschitz map, where $U\subset\R^k$ is open. Then for almost every $z\in U$ the metric directional derivative
 $\md\varphi_z(v)$ exists for every $v\in\R^k$. Furthermore, there are compact sets $K_i\subset U$, $i\in\N$, such that $\lm^k(U\ohne\cup K_i)=0$ and 
 such that the following property holds: 
 For every $i$ and every $\varepsilon>0$ there exists $r(K_i,\varepsilon)>0$ such that
 \begin{equation}\label{equation:strong-metric-norm}
  |d(\varphi(z+v),\varphi(z+w))-\md\varphi_z(v-w)|\leq \varepsilon|v-w|
 \end{equation}
 for all $z\in K_i$ and all $v,w\in\R^k$ satisfying $|v|,|w|\leq r(K_i,\varepsilon)$ and $z+w\in K_i$.
\et
Here $|\cdot|$ denotes the Euclidean norm and $\lm^k$ the Lebesgue measure.
If $\md\varphi_z(v)$ exists for all $v\in\R^k$ and satisfies \eqref{equation:strong-metric-norm} then $\md\varphi_z$  is called metric derivative of $\varphi$ 
at the point $z$. As a direct consequence of \thmref{theorem:metric-deriv-strong} the metric derivative is a seminorm. 
It is not difficult to prove that if $U\subset\R^k$ is merely Borel measurable then $\md\varphi_z$ can be defined at almost every Lebesgue density point $z\in U$ 
by a simple approximation argument.

\section{Proof of \thmref{theorem:main-thm-trees}}\label{section:charact-trees}
The following characterization of metric trees will be used in the proof of \thmref{theorem:main-thm-trees}.
\bp\label{proposition:trivial-loops}
 Let $X$ be a geodesic metric space. If for every Lipschitz curve $\gamma:[0,1]\to X$ with $\gamma(0)=\gamma(1)$ and all Lipschitz functions $f, \pi:X\to \R$ 
 \begin{equation}\label{equation:integral-loop-zero}
  \int_0^1 (f\circ\gamma)(t)\cdot(\pi\circ\gamma)'(t)\,dt= 0
 \end{equation}
 then $X$ is a metric tree.
\ep
Note that the function $\pi\circ\gamma$ is a Lipschitz function on $[0,1]$ and therefore, by Rademacher's theorem, its classical derivative 
$(\pi\circ\gamma)'(t)$ exists for almost every $t\in[0,1]$ and defines a bounded Borel function. In particular, the integral in \eqref{equation:integral-loop-zero}
is well-defined.

\begin{proof}
 We fist show that $X$ is uniquely geodesic. For this let $\gamma_1, \gamma_2:[0,1]\to X$ be two geodesics with the same endpoints. Set 
 $a:=d(\gamma_1(0), \gamma_1(1))$ and define a Lipschitz map
 by $\pi(x):= d(x, \gamma_1(0))$ and let $\gamma:[0,1]\to X$ be the constant speed loop given by the concatenation of $\gamma_1$ and 
 $\gamma_2^-(t):=\gamma_2(1-t)$. 
 Then $(\pi\circ\gamma)'(t)= 2a$ for every $t\in[0,\frac{1}{2}]$ and $(\pi\circ\gamma)'(t)= -2a$ for every $t\in[\frac{1}{2},1]$. Setting 
 $f(x):= \max\{0,1-\dist(\gamma_1([0,1]), x)\}$ it follows immediately from \eqref{equation:integral-loop-zero} that $\gamma_2([0,1])\subset\gamma_1([0,1])$ and therefore that $\gamma_2(t)= \gamma_1(t)$
 for all $t\in[0,1]$.
 
 We now prove that every geodesic triangle in $X$ is a tripod. For this let $\gamma_i:[0,a_i]\to X$, $i=1, 2, 3$, be three 
 geodesics parameterized by arc-length forming a geodesic triangle, i.e., such that $\gamma_1(a_1)=\gamma_2(0)$, 
 $\gamma_2(a_2)=\gamma_3(0)$, and $\gamma_3(a_3)=\gamma_1(0)$. 
 Set
 \begin{align*}
  t_1&:= \max\{t\in[0,a_1]: \gamma_1(t)\in \gamma_3([0,a_3])\}\\
  t_2&:= \min\{t\in[0,a_1]: \gamma_1(t)\in \gamma_2([0,a_2])\}.
 \end{align*}
 By defining $\gamma:[0,1]\to X$ to be the constant speed loop given by the concatenation of 
 $\gamma_1$, $\gamma_2$ and $\gamma_3$ and by setting $\pi(x):=d(x, \gamma_1(0))$ and
 \begin{equation*}
  f_\varepsilon(x):= \max\{0, 1-\varepsilon^{-1}\dist(x, \gamma_1([0,1]))\}
 \end{equation*}
 we easily find using the uniqueness of geodesics and the fact that
 \begin{equation*}
  \int_0^1 (f_\varepsilon\circ\gamma)(t)\cdot(\pi\circ\gamma)'(t)\,dt= 0
 \end{equation*}
 for every $\varepsilon>0$
 that $t_1=t_2$ and hence $\gamma_1(t_1)\in \image(\gamma_2)\cap\image(\gamma_3)$. By the uniqueness of geodesics the triangle consisting of 
 $\gamma_1, \gamma_2, \gamma_3$ is thus a tripod. This completes the proof.
\end{proof}

The proof of \thmref{theorem:main-thm-trees} involves several other auxiliary results which we state now. 
In order to simplify the language in what follows it is convenient to introduce a new terminology.
We say a metric space $X$ has property $(\Diamond)$ if it satisfies the hypotheses of \thmref{theorem:main-thm-trees}, that is, if 
there exists $\lambda\geq 1$ such that for any two closed balls $B_1, B_2\subset X$ with non-empty intersection 
there exist $z, z'\in X$ and $\nu\geq 0$ such that
 \begin{equation}\label{equation:intersection-property-before-special-case}
  B(z, \nu)\subset B_1\cap B_2\subset B(z', \lambda\nu).
 \end{equation}
To prove \thmref{theorem:main-thm-trees} it will in fact be enough to require \eqref{equation:intersection-property-before-special-case} only for those balls $B_1=B(x,r)$ and 
 $B_2=B(x',r')$ for which $\max\{r,r'\}<d(x,x')$. It is then readily verified that \eqref{equation:intersection-property-before-special-case} can be replaced by the hypothesis
 \begin{equation*}
  \diam (B(x,r)\cap B(x',r'))\leq 2\lambda(r+s-d(x,x')).
 \end{equation*}
We begin with the following observation.
\br\label{remark:diameter-estimate}{\rm
 Let $X$ be a geodesic metric space with property $(\Diamond)$. Then for any two distinct points $x,y\in X$ and for all $t\in (0,1)$ and 
 $0\leq h<\max\{t,1-t\}d(x,y)$ the set
\begin{equation*}
  A:= B(x, tr+h)\cap B(y,(1-t)r+h),
\end{equation*}
where $r:= d(x,y)$, satisfies $\diam A\leq 4\lambda h$.
}
\er
Choosing $h=0$ in the remark we obtain that such a space $X$ is uniquely geodesic.

\begin{proof}[{Proof of \remref{remark:diameter-estimate}}]
 By symmetry we may assume that $t\leq\frac{1}{2}$. Note that $A$ is non-empty since $X$ is geodesic; moreover $x,y\not\in A$.
 By assumption there exist $z, z'\in X$ and $s\geq 0$ such that $B(z, s)\subset A\subset B(z', s)$. It is clear that $s<d(x,z)$.
 Let $\alpha$ be a geodesic from $z$ to $x$ parameterized by arc-length. Then $\alpha(s)\in A$ and hence
 \begin{equation*}
  r\leq d(x,\alpha(s)) + d(\alpha(s), y)\leq tr+h-s+(1-t)r+h
 \end{equation*}
 from which it follows that $s\leq 2h$ and hence $\diam A\leq 2\lambda s\leq 4\lambda h$.
\end{proof}

\bl\label{lemma:bicombing-property}
 If $X$ is a geodesic metric space with property $(\Diamond)$ then for any points $x,y, y'\in X$ the  
 geodesics $c, c':[0,1]\to X$, parameterized proportional to arc-length, with $c(0)=x=c'(0)$ and $c(1)=y$ and $c'(1)=y'$ satisfy
 \begin{equation*}
  d(c(t), c'(t))\leq 4\lambda d(y,y')
 \end{equation*}
 for all $t\in[0,1]$.
\el
\begin{proof}
 Let $x,y,y'\in X$ be arbitrary points and set $r:=d(x,y)$ and $r':=d(x,y')$. If $y=x$ or $y'=x$ then the statement is trivially true and we may therefore 
 assume that $0<r'\leq r$. Set $h:= d(y,y')$ and fix $t\in(0,1)$. 
 If $h\geq tr$ then 
 \begin{equation*}
  d(c(t), c'(t))\leq d(c(t), x)+d(x, c'(t))= tr +tr'\leq 2h.
 \end{equation*}
 If, on the other hand, $h<tr$ then we define
 \begin{equation*}
  A:= B(x,tr+h)\cap B(y, (1-t)r+h)
 \end{equation*}
 and obtain that
 \begin{equation*}
  c'(t)\in B(x, tr')\cap B(y', (1-t)r') \subset B(x, tr)\cap B(y, (1-t)r+ h)\subset A
 \end{equation*}
 as well as $c(t)\in A$. Since by \remref{remark:diameter-estimate} we have  $\diam A\leq 4\lambda h$ it follows that
 \begin{equation*}
  d(c(t), c'(t))\leq 4\lambda d(y,y').
 \end{equation*}
 This concludes the proof.
\end{proof}
\br\label{remark:approximate-balls}{\rm\sloppy
 It should be noted that the assertions of \remref{remark:diameter-estimate} and \lemref{lemma:bicombing-property} remain true if $(\Diamond)$ is replaced by the 
 following condition: For any two balls $B(x, r), B(y, s)\subset X$ with non-empty intersection and every
 $\varepsilon>0$ there exist $z, z'\in X$ and $t\geq 0$ with
 \begin{equation*}
  B(z, t-\varepsilon)\subset B(x,r+\varepsilon)\cap B(y, s+\varepsilon)
 \end{equation*}
and
\begin{equation*}
 B(x, r)\cap B(y, s)\subset B(z', \lambda(t+\varepsilon)).
\end{equation*}
This property will be called $(\Diamond')$ and will be used in the proof of \corref{corollary:gromov-charact}.
}
\er
\fussy

As a consequence of the lemma above we obtain the following extension property.

\bc\label{corollary:lip-connected}
If $X$ is a geodesic metric space with property $(\Diamond')$ then every Lipschitz map $\varphi: S^m\to X$ has a Lipschitz extension 
$\overline{\varphi}: B^{m+1}\to X$ with $\lip(\overline{\varphi})\leq (8\lambda+12)\lip(\varphi)$.
In other words, $X$ is Lipschitz $m$-connected for every $m\geq 0$.
\ec

Here, $S^m$ and $B^{m+1}$ denote the unit sphere and unit ball of $\R^{m+1}$ with the Euclidean metric, respectively.

\begin{proof}
 For a given Lipschitz map $\varphi: S^m\to X$ fix $x_0\in \varphi(S^m)$ and define $\overline{\varphi}: B^{m+1}\to X$ by 
\begin{equation*}
 \overline{\varphi} (rz):= \left\{
 \begin{array}{l@{\quad:\quad}l}
  x_0 & 0\leq r\leq \frac{1}{2}\\
  c_{z}(2(r-\frac{1}{2})) & \frac{1}{2}<r\leq 1,
 \end{array}\right.
\end{equation*}
where $c_{z}:[0,1]\to X$ denotes the constant-speed geodesic from $x_0$ to $\varphi(z)$. This clearly defines a Lipschitz map which extends $\varphi$. A simple calculation
using \lemref{lemma:bicombing-property} shows furthermore that $\lip(\overline{\varphi})\leq (8\lambda+12)\lip(\varphi)$.
\end{proof}


%
In the proof of the following lemma we will use the obvious but useful fact that the properties $(\Diamond)$ and $(\Diamond')$ persist under rescaling of the metric.
\bl\label{lemma:bilipschitz-piece}
 Let $X$ be a geodesic metric space with property $(\Diamond')$ and let $K\subset \R^2$ be Borel measurable. 
 If $\varphi: K\to X$ is a Lipschitz map then $\md\varphi_z$ is degenerate for almost all $z\in K$.
\el
\begin{proof}
 We argue by contradiction and assume that $\md\varphi_z$ is non-degenerate, and thus a norm, for all $z$ in some set of positive Lebesgue 
 measure. By  \thmref{theorem:metric-deriv-strong} there exists a compact set $K'\subset K$ of positive Lebesgue measure and a Lebesgue density point $z_0\in K'$ 
 such that the norm $\|\cdot\|:=\md\varphi_{z_0}$ has the following property. For every $\varepsilon>0$ there exists $r_0>0$ such that the map 
 $\tilde{\varphi} : (K'-z_0, \|\cdot\|)\to X$ given by $\tilde{\varphi}(v):=\varphi(z_0+v)$ is $(1+\varepsilon)$-bilipschitz on $B(0, r_0)\cap (K'-z_0)$. 
 Here, $B(0, r_0)$ denotes the ball of radius $r_0$ with respect to the Euclidean metric whereas we will denote by $\hat{B}(0,r)$ the ball or radius $r$ with 
 respect to $\|\cdot\|$. Note that $\hat{B}(0,1)$ is convex, centrally symmetric and closed with respect to the
 Euclidean metric. Let $B(0,r_1)\subset \hat{B}(0,1)$ be the Euclidean ball of maximal radius and let $v_0\in\bdry B(0,r_1)\cap \bdry \hat{B}(0,1)$.
 Set $y:=2v_0$. It is then easy to see that
 \begin{equation}\label{equation:diam-ration-infinity}
 \frac{1}{h}\diam_{\|\cdot\|}\left[\hat{B}(0,1+h)\cap\hat{B}(y, 1+h)\right]\to\infty\qquad\text{as $h\to 0$.}
 \end{equation} 
 Indeed, this follows from the fact that
 \begin{equation*}
  \hat{B}(0,1+h)\cap\hat{B}(y, 1+h)\supset B(0,r_1(1+h))\cap B(y, r_1(1+h))
 \end{equation*}
 and
 \begin{equation*}
 \diam_{\text{eucl.}}\left[B(0,r_1(1+h))\cap B(y, r_1(1+h))\right] = 2\sqrt{2h+h^2}r_1.
 \end{equation*}
 Since $z_0$ is a Lebesgue density point of $K'$ there exist sequences $r_n\searrow 0$ and $y_n\in\R^2$ such that $y_n\in K_n:=\frac{1}{r_n}(K'-z_0)$ and 
 $s_n:=\|y_n-y\|\to 0$ and such that for every $h>0$ 
  \begin{equation}\label{equation:Hausdorff-convergence}
   \left[\hat{B}(0,1+h)\cap\hat{B}(y_n, 1+h)\cap K_n\right]\longrightarrow  \left[\hat{B}(0,1+h)\cap\hat{B}(y, 1+h)\right]
 \end{equation}
 with respect to Hausdorff convergence.
 Denote by $X_n$ the metric space $(X, r_n^{-1}d)$ and note that $X_n$ has property $(\Diamond')$. Then the maps
 $\tilde{\varphi}_n:(K_n, \|\cdot\|)\to X_n$ given by $\tilde{\varphi}_n(u):= \tilde{\varphi}(r_n u)$ are $(1+\varepsilon_n)$-bilipschitz with 
 $\varepsilon_n\searrow0$, see \thmref{theorem:metric-deriv-strong}. 
 \remref{remark:diameter-estimate} then implies that for $h>0$ small enough
 \begin{equation}\label{equation:bounded-diam-ratio}
  \diam_{X_n}\left[B_{X_n}\left(\tilde{\varphi}_n(0), t_n+ h\right)\cap B_{X_n}(\tilde{\varphi}_n(y_n), t_n+h)\right]\leq 4\lambda h
 \end{equation}
 for every $n\in\N$, where $t_n:= \frac{1}{2}d_{X_n}(\tilde{\varphi}_n(0), \tilde{\varphi}_n(y_n))$.
 On the other hand, it follows from \eqref{equation:Hausdorff-convergence} that for $n$ large enough
 \begin{equation*}
  \begin{split}
  \diam&_{X_n}\left[B_{X_n}\left(\tilde{\varphi}_n(0), t_n+h\right)\cap B_{X_n}(\tilde{\varphi}_n(y_n), t_n+h)\right]\\
   &\geq \frac{1}{1+\varepsilon_n}\diam_{\|\cdot\|}\left[\hat{B}\left(0, (1+\varepsilon_n)^{-1}(t_n+h)\right)\cap \hat{B}\left(y_n, (1+\varepsilon_n)^{-1}(t_n+h)\right)\cap K_n\right]\\
   &\geq  \frac{1}{2(1+\varepsilon_n)}\diam_{\|\cdot\|}\left[\hat{B}\left(0, 1+h'_n\right)\cap \hat{B}\left(y_n,1+h'_n\right)\right]\\ 
 \end{split}
 \end{equation*}
 where
 \begin{equation*}
  h'_n=\frac{h-\frac{s_n}{2} - 2\varepsilon_n - \varepsilon_n^2}{(1+\varepsilon_n)^2}.
 \end{equation*}
 Using the fact that $\varepsilon_n, s_n\to 0^+$ together with \eqref{equation:diam-ration-infinity} we readily arrive at a contradiction with 
 \eqref{equation:bounded-diam-ratio}.
 This completes the proof.
 \end{proof}

We are finally ready to prove \thmref{theorem:main-thm-trees}.
\begin{proof}[{Proof of \thmref{theorem:main-thm-trees}}]
 We first show that if $X$ has property $(\Diamond')$ then $X$ is a metric tree. We do this by verifying the hypothesis \eqref{equation:integral-loop-zero}
 of \propref{proposition:trivial-loops}. For this let $\gamma: S^1\to X$ be a Lipschitz loop and let $\tilde{\gamma}: B^2\to X$ be 
 a Lipschitz extension, which exists by \corref{corollary:lip-connected}. 
 Let $f, \pi: X\to \R$ be arbitrary Lipschitz functions. By Stokes' Theorem and an obvious smoothing argument we obtain
 \begin{equation*}
  \int_{S^1}(f\circ\gamma)(t)\cdot(\pi\circ\gamma)'(t)dt =\int_{B^2}\det(D_z[(f,\pi)\circ\tilde{\gamma}]dz= 0, 
 \end{equation*}
where the second equality follows from \lemref{lemma:bilipschitz-piece} and \remref{remark:det-degenerate}. It now follows from \propref{proposition:trivial-loops} that $X$ is a metric tree.
The proof of the converse statement is a consequence of the lemma below with $\delta=0$.
 \end{proof}

\bl\label{Lemma:hyp-distortion}
 Let $\delta\geq 0$ and suppose $(X,d)$ is a geodesic metric space all of whose geodesic triangles are $\delta$-thin. 
 Then for any two closed balls $B_1,B_2$ in $X$ with non-empty intersection there exist $z\in X$ and $\nu\geq 0$ such that
 \begin{equation*}
  B(z, \nu)\subset B_1\cap B_2\subset B(z,\nu+\delta).
 \end{equation*}
\el

The proof of the lemma is straight-forward and is implicitly contained in \cite{Chatterji-Niblo}. For completeness we give our own short proof here.

\begin{proof}
 Let $B(x,r)$ and $B(y,s)$ be two closed balls in $X$ and suppose $K:= B(x,r)\cap B(y,s)$ is non-empty. 
 We may assume without loss of generality that $r\geq s$. Set 
 $d:=d(x,y)$. If $r-s>d$ then it follows that $K=B(y,s)$ and there is nothing to prove.
 If $r-s\leq d$ then choose a geodesic $c$ between $x$ and $y$ and let $z\in\image(c)$ be such that
 \begin{equation*}
  d(x,z) = \frac{r-s+d}{2}\quad\text{and}\quad d(y,z)=\frac{s-r+d}{2}.
 \end{equation*}
 Set furthermore $\nu:= \frac{1}{2}(r+s-d)$. Since $K\not=\emptyset$ we clearly have $\nu\geq 0$. We claim that 
 \begin{equation}\label{equation:distortion-inclusion}
  B(z,\nu)\subset K\subset B(z, \nu+\delta).
 \end{equation}
 To prove the first inclusion it is enough to note that for $w\in B(z, \nu)$ 
 \begin{equation*}
  d(w,x)\leq d(w,z)+d(z,x)\leq \nu+\frac{r-s+d}{2}= r
 \end{equation*}
 and analogously
 \begin{equation*}
  d(w,y)\leq d(w,z)+d(z,y)\leq \nu+\frac{s-r+d}{2}= s.
 \end{equation*}
 As for the proof of the second inclusion we let $w\in K$ and consider a geodesic triangle with vertices $x_1:=x$, $x_2=y$, $x_3=w$ and 
 edges $c_{k\ell}$, where $c_{12}:= c$. Let $a_1, a_2, a_3$ be as in the definition of the $\delta$-thinness. If $d(x,z)\leq a_1$ 
 then
 \begin{equation*}
  d(w, z)\leq d(w,c_{13}(d(x,z)))+\delta= d(w,x)-d(x,z)+\delta\leq r-\frac{r-s+d}{2+\delta}=\nu+\delta,
 \end{equation*}
 whereas if $d(x,z)>a_1$ we compute
 \begin{equation*}
  d(w,z)\leq d(w, c_{23}(d(y,z)))+\delta=d(w,y)-d(y,z)+\delta\leq s-\frac{s-r+d}{2}+\delta=\nu+\delta.
 \end{equation*}
 This establishes the second inclusion in \eqref{equation:distortion-inclusion} and completes the proof.
\end{proof}
%
%
%
%
%
%

\section{Proof of \corref{corollary:gromov-charact}}
The following proof uses asymptotic cones. For definitions and basic properties we refer to \cite{Bridson-Haefliger}. We will need the following crucial fact:
A geodesic metric space $X$ is Gromov hyperbolic if and only if every asymptotic cone of $X$ is a metric tree, see e.g.~\cite{Drutu}.

\begin{proof}[{Proof of \corref{corollary:gromov-charact}}]
 By Proposition 3.A.1 in \cite{Drutu} it suffices to show that every asymptotic cone of $X$ is a metric tree.  
 Let therefore $\omega$ be a non-principal ultrafilter on $\N$, $(x_n)\subset X$ and $r_n\nearrow\infty$ and denote by $X_\omega$ the asymptotic cone associated to 
 the sequence $(X, \frac{1}{r_n}d, x_n)$ and $\omega$.
 Clearly, $X_\omega$ is geodesic; thus we only need to show that $X_\omega$ has property $(\Diamond')$.
 For this let $\overline{y}=(y_n)$ and $\overline{z}=(z_n)$ be arbitrary points in $X_\omega$ and let $r, s>0$ be such that $r+s\geq d_\omega(\overline{y}, \overline{z})$. 
 Fix $\varepsilon >0$. By assumption there exist for $\omega$-almost every $n\in\N$ two points $u_n, u'_n\in X$ and $t_n\geq 0$ such that
 \begin{equation}\label{equation:intersection-property-downstairs}
  B(u_n, t_n)\subset B(y_n, (r+\varepsilon)r_n)\cap B(z_n, (s+\varepsilon)r_n)\subset B(u'_n, \lambda t_n+\delta).
 \end{equation}
 We have $t_n\leq (r+\varepsilon)r_n$ from which we conclude that $\bar{t}:=\lim_\omega\frac{t_n}{r_n}$ exists.
 By the triangle inequality we furthermore have that $\overline{u}:=(u_n)$ and $\overline{u}':=(u'_n)$ are elements of $X_\omega$. We first claim that
 \begin{equation*}
  \hat{B}(\overline{u}, \bar{t}-\varepsilon)\subset \hat{B}(\overline{y}, r+\varepsilon)\cap \hat{B}(\overline{z}, s+\varepsilon)
 \end{equation*}
 for every $\varepsilon>0$ small enough. Here, $\hat{B}$ denotes a ball in $X_\omega$.
 Let $\overline{v}= (v_n)\in \hat{B}(\overline{u}, \bar{t}-\varepsilon)$ and note that for $\omega$-almost all $n\in\N$ we have $d(v_n, u_n)\leq t_n$ and therefore,
 by \eqref{equation:intersection-property-downstairs}
 \begin{equation*}
  \frac{1}{r_n}d(v_n, y_n)\leq r+\varepsilon\qquad\text{and}\qquad \frac{1}{r_n}d(v_n, z_n)\leq s+\varepsilon
 \end{equation*}
 for $\omega$-almost all $n\in\N$.
 This clearly implies that $\overline{v}\in \hat{B}(\overline{y}, r+\varepsilon)\cap \hat{B}(\overline{z}, s+\varepsilon)$.
 Next, we show that 
  \begin{equation}\label{equation:second-inclusion}
   \hat{B}(\overline{y}, r)\cap \hat{B}(\overline{z}, s)\subset \hat{B}(\overline{u}', \lambda\bar{t}).
 \end{equation}
 Indeed, if $\overline{v}\in B(\overline{y}, r)\cap B(\overline{z}, s)$ then $v_n\in B(y_n, (r+\varepsilon)r_n)\cap B(z_n, (s+\varepsilon)r_n)$ for 
 $\omega$-almost all $n\in\N$ and therefore
 \begin{equation*}
  d_\omega(\overline{v}, \overline{u}')= \lim_\omega \frac{1}{r_n}d(v_n, u'_n)\leq \lim_\omega\frac{\lambda t_n+\delta}{r_n}= \lambda \bar{t},
 \end{equation*}
 which proves \eqref{equation:second-inclusion}.
 By \thmref{theorem:main-thm-trees}, which, as we proved, holds when $(\Diamond)$ is replaced by $(\Diamond')$, this is enough to conclude that $X_\omega$ is a metric tree. 
\end{proof}

\section{Appendix: A proof of \thmref{theorem:metric-deriv-strong}}

A simple proof of the following fact can be found e.g.~in \cite[Theorem 4.1.6]{Ambrosio-Tilli}.

\bp\label{proposition:one-dimensional-derivative}
 Let $X$ be a metric space and $\gamma: (a,b)\to X$ a Lipschitz curve. Then $ \md\varphi_t(1)$ exists for almost every $t\in(a,b)$.
\ep

\begin{proof}[Proof of \thmref{theorem:metric-deriv-strong}]
 Observe first that if $z\in U$ and $r>0$ and if $v, v'\in\R^k$ are such that $z+rv, z+rv'\in U$ then
 \begin{equation}\label{equation:double-ineq-triang}
  \frac{1}{r}d(\varphi(z+rv'),\varphi(z)) -\varrho\leq \frac{1}{r}d(\varphi(z+rv),\varphi(z))\leq\frac{1}{r}d(\varphi(z+rv'),\varphi(z)) +\varrho
 \end{equation}
 where $\varrho:=\lip(\varphi)|v'-v|$.
Moreover, if $z\in U$ and $v,v'\in\R^k$ are such that $\md\varphi_z(v)$ and $\md\varphi_z(v')$ exist then
 \begin{equation}\label{equation:lip-md}
  |\md\varphi_z(v)-\md\varphi_z(v')|\leq\lip(\varphi)|v-v'|. 
 \end{equation}

We can now easily show that for $\lm^k$-almost every $z\in U$ the limit $\md\varphi_z(v)$ exists for all $v\in\R^k$. Indeed,
 let $(v_i)_{i\in\N}\subset S^{k-1}$ be a 
 countable
 dense subset. By \propref{proposition:one-dimensional-derivative} and Fubini's theorem there exists an $\lm^k$-negligible set $N\subset U$ such 
 that $\md\varphi_z(v_i)$ exists for all $i\in\N$ and all $z\in U\ohne N$.
 Given $v\in S^{k-1}$ arbitrary choose a subsequence $v_{i_j}$ of $(v_i)$ which converges to $v$.  The existence of $\md\varphi_z(v)$ now follows immediately from
 \eqref{equation:double-ineq-triang} and \eqref{equation:lip-md}. Moreover, $\md\varphi_z(sv)$ exists for all 
 $s\geq 0$ and $\md\varphi_z(sv)=s\md\varphi_z(v)$. Hence, $\md\varphi_z(v)$ exists for all $v\in\R^k$ and all $z\in U\ohne N$.

We turn to the proof of \eqref{equation:strong-metric-norm}.
Denote by $S^{k-1}$ the unit sphere in $\R^k$ and by $C(S^{k-1})$ the space of continuous real-valued functions on $S^{k-1}$, endowed with the supremum norm.
Denote by $D$ the set of points $z\in U$ where $\md\varphi_z(v)$ exists for all $v\in\R^k$.
If $z\in D$ then, by  \eqref{equation:lip-md}, the function $f_z:S^{k-1}\to[0,\infty)$ given by $f_z(v):=\md\varphi_z(v)$ is $\lip(\varphi)$-Lipschitz.  
Furthermore, the map $F: D\to C(S^{k-1})$ given by $F(z):= f_z$ is $\lm^k$-measurable. Since
$C(S^{k-1})$ is separable, 
 Lusin's Theorem \cite[2.3.5]{Federer} asserts the existence of compact subsets $K'_i\subset D$ with $\lm^k(D\ohne\cup K'_i)=0$ and such that $F$ is continuous on 
 each $K'_i$. 
 In particular,  the function $g_r(z):= \sup_{|v|=1}|\frac{1}{r}d(\varphi(z+rv),\varphi(z))-\md\varphi_z(v)|$ is continuous on each $K'_i$, and
 converges pointwise
 to $0$ as $r\searrow0$. By Egoroff's Theorem \cite[2.3.7]{Federer} there exist compact subsets $K_i^j\subset K'_i$ with 
 $\lm^k(D\ohne\cup_{i,j} K_i^j)=0$ and such that $g_r$ converges to $0$ uniformly on each $K_i^j$ as $r\searrow 0$. In particular, for every $\varepsilon>0$ there exists
 $r(K_i^j,\varepsilon)>0$ such that $g_r(z)\leq \varepsilon/2$ and $\|\md\varphi_z-\md\varphi_{z'}\|_{\infty}\leq \varepsilon/2$ whenever $0<r\leq 2r(K_i^j,\varepsilon)$ and
 $z, z'\in K_i^j$ with $|z-z'|\leq r(K_i^j,\varepsilon)$. If $z\in K_i^j$ and if $v,w\in\R^k$  satisfy $|v|,|w|\leq r(K_i^j,\varepsilon)$, $v\not= w$, and $z+w\in K_i^j$
 then, setting $r:=|v-w|$, $z':= z+w$, and $v':= \frac{1}{r}(v-w)$, we obtain
 \begin{equation*}
  \begin{split}
   |d(\varphi(z+v),&\varphi(z+w))-\md\varphi_z(v-w)|\\ 
                                                    &\leq r\left|\frac{1}{r}d(\varphi(z'+rv'),\varphi(z'))-\md\varphi_{z'}(v')\right| +\frac{\varepsilon}{2}|v-w|\\
                                                    &\leq r g_r(z')+\frac{\varepsilon}{2}|v-w|\\
                                                    &\leq \varepsilon|v-w|.\\
  \end{split}
 \end{equation*}
 This proves \eqref{equation:strong-metric-norm} and hence the theorem.
\end{proof}


\begin{thebibliography}{331}
  \bibitem{Ambrosio-Tilli}L.~Ambrosio, P.~Tilli: {\it Topics on Analysis in metric spaces}, Oxford Lecture Series in Mathematics and its Applications, 25. Oxford University Press, Oxford, 2004.
  \bibitem{Bridson-Haefliger}M.~R.~Bridson, A.~Haefliger: {\it Metric Spaces of Non-Positive Curvature}, 
    Grundlehren der mathematischen Wissenschaften 319, Springer, 1999.
  \bibitem{Chatterji-Niblo}I.~Chatterji, G.~A.~Niblo: {\it A characterization of hyperbolic spaces}, Groups, Geometry, and Dynamics, 1 (2007), no. 3, 281Ð-299.
    \bibitem{Drutu}C.~Dru\c{t}u: {\it Quasi-isometry invariants and asymptotic cones}, International Conference on Geometric and Combinatorial Methods in Group Theory 
   and Semigroup Theory (Lincoln, NE, 2000).  Internat. J. Algebra Comput.  12  (2002),  no. 1-2, 99--135. 
  \bibitem{Federer}H.~Federer: {\it Geometric measure theory}, Die Grundlehren der mathematischen Wissenschaften, Band 153 Springer-Verlag New York Inc., New York 1969 and 1996.
  \bibitem{Ghys-delaHarpe}E.~Ghys, P.~de la Harpe (ed): {\it Sur les groupes hyperboliques d'apr\`es Mikhael Gromov}, Progr. Math. vol 83, 
    Birkh\"auser, Boston MA, 1990.
  \bibitem{Gromov-hyperbolic}M.~Gromov: {\it Hyperbolic groups}, Essays in group theory,  75--263, Math. Sci. Res. Inst. Publ., 8, Springer, New York, 1987.
  \bibitem{Preiss-Schechtman}W.~B.~Johnson, J.~Lindenstrauss, D.~Preiss, G.~Schechtman: {\it Lipschitz quotients from metric trees and Banach spaces containing $\ell_1$}, 
   J. Funct. Anal. 194, (2002), no. 2, 332 -- 346.
  \bibitem{Kirchheim} B.~Kirchheim: {\it Rectifiable metric spaces: local structure and regularity of the Hausdorff measure}, 
    Proc. Am. Math. Soc. 121 (1994), no. 1, 113--123.
  \bibitem{Korevaar-Schoen}N.~J.~Korevaar, R.~M.~Schoen: {\it Sobolev spaces and harmonic maps for metric space targets}, 
   Comm. Anal. Geom. 1 (1993), no. 3-4, 561--659.
  \bibitem{Wenger-sharp-constant}S.~Wenger: {\it Gromov hyperbolic spaces and the sharp isoperimetric constant}, preprint 2007.
 \end{thebibliography}
\end{document}